\documentclass[11pt,a4paper]{article}
\usepackage{amssymb}
\usepackage{amsmath}
\usepackage{amsthm}
\theoremstyle{plain}
\newtheorem{thm}{Theorem}[section]

\newtheorem{lem}[thm]{Lemma}

\newtheorem{prop}[thm]{Proposition}
 
\newtheorem{rem}{Remark}[section]

\def\bproof{\textbf{Proof}: }
\def\eproof{\hfill$\Box$}
\def \F{\mathbb F}

\author{ Ay\c ca \c Ce\c smelio\u glu \\
\small Sabanc\i \; University \\
\small MDBF, Orhanl\i, Tuzla \\
\small   34956 \.{I}stanbul, Turkey \\
\small cesmelioglu@sabanciuniv.edu}

\title{A Representation of Permutations with Full Cycle 
}
\date{ }

\begin{document}

\maketitle

\begin{abstract}
For $q>2$, Carlitz proved that the group of permutation polynomials (PPs) over $\F_q$ is generated by linear polynomials and $x^{q-2}$. 
Based on this result, this note points out a simple method for representing all PPs with full cycle 
over the prime field $\F_p$, where $p$ is an odd prime. We use the isomorphism between the symmetric group $S_p$ of $p$ elements
and the group of PPs over $\mathbb{F}_p$, and the well-known fact that permutations in $S_p$ have the same cycle structure if and only if 
they are conjugate.
\end{abstract}

\section{Introduction}
For pseudorandom number generation, one of the important problems is to generate sequences with long periods.
A sequence $(s_n)$, which is recursively defined by the formula $s_{n+1}=\psi(s_n)$ with an initial value $s_0$,
will be a purely periodic sequence when $\psi$ is a PP over the finite field $\F_q$ with 
$q=p^r$ elements where $p$ is a prime and $r\geq 1$. The period length of $(s_n)$ equals the length of the 
cycle in which the initial value $s_0$ lies in the permutation induced by $\psi$.  
The maximal possible period $q$ for $(s_n)$ will be achieved, if $\psi$ corresponds to a PP with full cycle.
Note that we use the term \textit{PP with full cycle} for a PP which induces a full cycle permutation.

For $q \geq 2$, PPs over $\F_q$ form a group under composition and reduction modulo $x^q-x$ 
and this group is isomorphic to the symmetric group $S_q$.
In \cite{carlitz}, Carlitz proved that any transposition $(0a)$, for $a\in \F_q^*, q>2$, can be represented by the polynomial \[-a^2(((x-a)^{q-2}+a^{-1})^{q-2}-a)^{q-2},\] 
which shows that any PP over $\F_q$ is of the form
\begin{equation}\label{Pn}
\mathcal{P}_n(x)=(\cdots((a_0x+a_1)^{q-2}+a_2)^{q-2}+\cdots + a_n){q-2}+a_{n+1}
\end{equation}
for some $n \geq 1$ and $a_0,a_1,\cdots,a_{n+1}\in\F_q$ with $a_0\neq 0$ or a linear polynomial $\mathcal{P}_0(x)=cx+d \in \F_q[x]$.
The cycle structure of PPs given as in (\ref{Pn}) was studied in \cite{aaw} and also some conditions for obtaining PPs with full cycle were determined.

In this note, we characterize all PPs over $\F_p$ with full cycle, where $p$ is an odd prime, in terms of the representation as in (\ref{Pn}). As far as we know, no complete characterization of full cycle permutations was given in terms of PPs before.

For the following parts, we will denote the cycle decomposition of a permutation $\mathcal{P}$,
which can be expressed as a product of disjoint cycles as 
\[\mathcal{P}=\tau_1^{(1)} \tau_2^{(1)} \cdots \tau_{n_1}^{(1)} \tau_1^{(2)} \tau_2^{(2)} \cdots \tau_{n_2}^{(2)} \cdots \tau_{1}^{(s)} \tau_2^{(s)} \cdots \tau_{n_s}^{(s)}\]
where $\tau_{j}^{(i)}$, $ 1 \leq j \leq n_i$, is a cycle of length $\ell_i$, by
\[\mathcal{T}(\mathcal{P}) = [n_1\times \ell_1,n_2\times \ell_2, . . . , n_s\times \ell_s].\]

\section{Main Result}

The result of this note is based on the cycle decomposition of the permutations induced by linear polynomials and the following well-known proposition, for which we refer to Section 4.3 of \cite{dummit-foote}.
\begin{prop}\label{conjugacy}
Two elements of $S_q$ are conjugate in $S_q$ if and only if they have the same cycle structure.
\end{prop}

Note that for the linear polynomial $\mathcal{P}_0(x)=cx+d \in \F_q[x]$, the cycle decomposition is given by
\begin{equation}\label{LP-cycle decomposition}
\mathcal{T}(\mathcal{P}_0)=\left\{\begin{array}{ll}
                           \left[\frac{q}{p} \times p \right] & \;\text{if}\; c=1 \;\text{and}\; d\neq 0, \\
                           & \\
                           \left[\frac{q-1}{k} \times k, 1\times 1\right] &\; \text{if}\; c\neq 1,
                           \end{array}\right.
\end{equation}
where $k$ denotes the order of $c$ in $\F_q^*$.

In order to use the previous proposition, we need a representation for the PP which induces the inverse of a given
permutation. For the linear polynomial $\mathcal{P}_0(x)$, the inverse is given by $\mathcal{P}_0^{-1}(x)=c^{-1}x-c^{-1}d$.

In the following parts, we assume $q=p^r, r\geq 1$, and $q > 2$.

\begin{lem}\label{inverse}
Let $\mathcal{P}_n(x) \in \F_q[x]$ denote the PP
\[\mathcal{P}_n(x)=(\cdots((a_0x+a_1)^{q-2}+a_2)^{q-2}+\cdots+a_n)^{q-2}+a_{n+1},\]
where $a_0\in \F_q^*, a_1, a_2,\cdots,a_{n+1}\in \F_q, n\geq 1.$\\
Then the coefficients of 
\[\mathcal{P}_n^{-1}(x)=(\cdots((b_0x+b_1)^{q-2}+b_2)^{q-2}+\cdots+b_n)^{q-2}+b_{n+1}\] 
are given by
\[b_k=\left\{\begin{array}{ll}-a_0^{-1}a_{n+2-k}& \text{if n+2-k is odd},\\ 
-a_0a_{n+2-k} & \text{if n+2-k is even},\end{array}\right.\]
for $1\leq k \leq n+1$ and 
\[b_0=\left\{\begin{array}{ll} a_0 & \text{if n is odd},\\
                          a_0^{-1} & \text{if n is even}. 
             \end{array}\right.\]
\end{lem}

\bproof

Let 
\[\mathcal{P}_n^{-1}(x)=(\cdots((b_0x+b_1)^{q-2}+b_2)^{q-2}+\cdots+b_n)^{q-2}+b_{n+1}.\]
Then
\[(\mathcal{P}_n \circ \mathcal{P}_n^{-1})(x)=(\cdots((a_0\mathcal{P}_n^{-1}(x)+a_1)^{q-2}+a_2)^{q-2}+\cdots+a_n)^{q-2}+a_{n+1}\]
\[=(\cdots((a_0(\cdots((b_0x+b_1)^{q-2}+b_2)^{q-2}+\cdots+b_n)^{q-2}+a_0b_{n+1}+a_1)^{q-2}+a_2)^{q-2}+\cdots+a_n)^{q-2}+a_{n+1}.\]
We can write $(\mathcal{P}_n \circ \mathcal{P}_n^{-1})(x)$ as
\[(\cdots(((\cdots((a_0b_0x+a_0b_1)^{q-2}+a_0^{-1}b_2)^{q-2}+\cdots+a_0^{-1}b_n)^{q-2}+a_0b_{n+1}+a_1)^{q-2}+a_2)^{q-2}+\cdots+a_n)^{q-2}
 +a_{n+1}\]
if $n$ is even, and as
\[(\cdots(((\cdots((a_0^{-1}b_0x+a_0^{-1}b_1)^{q-2}+a_0b_2)^{q-2}+\cdots+a_0^{-1}b_n)^{q-2}+a_0b_{n+1}+a_1)^{q-2}+a_2)^{q-2}+\cdots+a_n)^{q-2}
 +a_{n+1}\]
when $n$ is odd.\\
We only consider the case where $n$ is odd since the case when $n$ is even can be proved similarly.
If we set 
\[b_{n+1}=-a_0^{-1}a_1,\; b_n=-a_0a_2,\; b_{n-1}=-a_0^{-1}a_3,\cdots,b_2=-a_0^{-1}a_n,\; b_1=-a_0a_{n+1},\;b_0=a_0\] 
then it is easy to see that 
\[(\mathcal{P}_n \circ \mathcal{P}_n^{-1})(x)=x\]
by inserting the coefficients $b_i, i=0,1,\cdots,n+1$, to the previous equation.

\eproof

The following theorem gives a complete description of the PPs over $\F_p$ with full cycle.
\begin{thm}

Let $p$ be an odd prime, $\widetilde{\mathcal{P}}_0(x)$ denote the linear polynomial of the form $x+d\in\F_p[x], d \in \F_p^*$, and $\widetilde{\mathcal{P}}_{2n}(x),n\geq 1$, denote a PP of the form
\[\widetilde{\mathcal{P}}_{2n}(x)=(\cdots(((\cdots((x+a_1)^{p-2}+a_2)^{p-2}+\cdots+a_n)^{p-2}+a_{n+1})^{p-2}-a_{n})^{p-2}-\cdots-a_2)^{p-2}-a_1\]
where $a_1,a_2,a_3,\cdots,a_{n+1}\in\F_p, a_{n+1}\neq 0$. 
The PP $\mathcal{P}(x)\in\F_p[x]$ is a full cycle if and only if it has a representation as $\mathcal{P}(x)=\widetilde{\mathcal{P}}_{2n}(x)$ for some $n \geq 0$.
\end{thm}
\bproof

The linear polynomial $x+d\in\F_p[x]$ with $d\neq 0$ corresponds to a permutation with full cycle and one can 
write any PP with full cycle as a conjugate $(\mathcal{P}_n\circ(x+d)\circ\mathcal{P}_n^{-1})(x)$ 
by using some $\mathcal{P}_n(x)\in \mathbb{F}_p[x], n\geq 0$.\\
Here, we only consider the case when $n\geq 1$ is even, the case where $n$ is odd can also be dealt with similarly and the case $n=0$
gives rise to linear PPs which are a full cycle.\\
Let \[\mathcal{P}_n(x)=(\cdots((b_0x+b_1)^{p-2}+ b_2)^{p-2}+\cdots+ b_n)^{p-2}+ b_{n+1}\in\mathbb{F}_p[x].\]
If $n\geq 2$ is even, then we have
\[\mathcal{P}_n^{-1}(x)=(\cdots((b_0^{-1}x-b_0^{-1}b_{n+1})^{p-2}-b_0b_n)^{p-2}-\cdots-b_0b_2)^{p-2}-b_0^{-1}b_1.\]
Therefore 
\[(\mathcal{P}_n\circ(x+d)\circ\mathcal{P}_n^{-1})(x)\]
\[=(\cdots((b_0(\cdots((b_0^{-1}x-b_0^{-1}b_{n+1})^{p-2}-b_0b_n)^{p-2}-\cdots-b_0b_2)^{p-2}-b_1+b_0d+b_1)^{p-2}+b_2)^{p-2}\]
\[+\cdots+b_n)^{p-2}+b_{n+1}\]
\[=(\cdots(((\cdots((x-b_{n+1})^{p-2}-b_n)^{p-2}-\cdots-b_2)^{p-2}+b_0d)^{p-2}+b_2)^{p-2}+\cdots+b_n)^{p-2}+b_{n+1}.\]
Note that the coefficients of the resulting PP satisfy the conditions of the theorem.

Conversely, for any given PP of the form $\widetilde{\mathcal{P}}_{2n}(x)$, we can find the linear polynomial $x+d$ and the PP $\mathcal{P}_n(x)$
such that $\widetilde{\mathcal{P}}_{2n}(x)=\mathcal{P}_n(x)\circ (x+d)\circ \mathcal{P}_n^{-1}(x)$ by tracing back the previous part of the proof. 
\eproof

\begin{rem} 
The kth iterate of $\widetilde{\mathcal{P}}_{2n}(x)\in \F_p[x]$ is easily seen to be
\[\widetilde{\mathcal{P}}_{2n}^{(k)}(x)=(\cdots(((\cdots((x+a_1)^{p-2}+a_2)^{p-2}+\cdots+a_n)^{p-2}+ka_{n+1})^{p-2}-a_{n})^{p-2}-\cdots-a_2)^{p-2}-a_1,\]
giving the elements of the Sylow p-subgroup generated by $\widetilde{\mathcal{P}}_{2n}(x)$.
\end{rem}

In terms of applications, PPs with full cycle are certainly the most interesting case.
However by using Proposition \ref{conjugacy}, and the isomorphism between the group of PPs over $\F_q$ and the symmetric group $S_q$, it becomes easy to represent permutations in $S_q$ having the same cycle structure as one of the linear polynomials.
We give the representation of such polynomials without a proof, since the calculations are easy but rather long.

\begin{rem}
Let $\widetilde{\mathcal{P}}_{2n}(x) \in \F_q[x]$ be the PP
\[\widetilde{\mathcal{P}}_{2n}(x)=(\cdots((b_0x+b_1)^{q-2}+b_2)^{q-2}+\cdots+b_{2n})^{q-2}+b_{2n+1},\]
with the coefficients 
\begin{equation}\label{Fqrepr}
b_i=\left\{\begin{array}{ll}
              c           & \text{for}\; i=0,\\
              ca_{i}      & \text{for}\; 1 \leq i \leq n\; \text{and i is odd},\\
              c^{-1}a_{i} & \text{for}\; 1 \leq i \leq n\; \text{and i is even},\\
              a_{n+1}     & \text{for}\; i=n+1,\\
              -a_{2n+2-i} & \text{for}\; n+2 \leq i \leq 2n+1,
             \end{array}\right.
\end{equation}
for some $c \in \F_q^*$ and $a_1,a_2,\cdots,a_{n+1} \in \F_q$.\\
By Proposition \ref{conjugacy} and Lemma \ref{inverse}, any permutation of $S_q$ having the same cycle decomposition as the linear polynomial $\mathcal{P}_0(x)=cx+d$ can be represented as in (\ref{Fqrepr}) for some $n \geq 0.$\\
Furthermore, the kth iterate of $\widetilde{\mathcal{P}}_{2n}(x)$ is
\[\widetilde{\mathcal{P}}_{2n}^{(k)}(x)=(\cdots((\beta_0x+\beta_1)^{q-2}+\beta_2)^{q-2}+\cdots+\beta_{2n})^{q-2}+\beta_{2n+1}\]
with 
\[\beta_i=\left\{\begin{array}{ll}
              c^k          & \text{for}\; i=0,\\
              c^ka_{i}      & \text{for}\; 1 \leq i \leq n\; \text{and i is odd},\\
              c^{-k}a_{i} & \text{for}\; 1 \leq i \leq n\; \text{and i is even},\\
              (1+c^{-1}+\cdots+c^{-(k-1)})a_{n+1}     & \text{for}\; i=n+1,\\
              -a_{2n+2-i} & \text{for}\; n+2 \leq i \leq 2n+1,
             \end{array}\right.\]
If $c=1$ and $\widetilde{\mathcal{P}}_{2n}(x)\neq x$ then $\mathcal{T}(\widetilde{\mathcal{P}}_{2n})=\left[\frac{q}{p} \times p \right]$, and if $c\neq 1$, then $\mathcal{T}(\widetilde{\mathcal{P}}_{2n})=\left[\frac{q-1}{k} \times k, 1\times 1\right]$, where $k=ord(c)$, by (\ref{LP-cycle decomposition}).

\end{rem}

\section{Acknowledgement}
The author would like to thank Wilfried Meidl for his careful reading and useful comments.

\end{document}